\input amstex
\documentstyle{amsppt}
\pagewidth{6.5in}
\pageheight{9in}
\topmatter
\title
$L^1$ is complemented in the dual space $L^{\infty *}$
\endtitle
\author Javier Guachalla H.\endauthor
\affil Universidad Mayor de San Andr\'es\footnote{e-mail:jguachallah\@gmail.com} \\La Paz - Bolivia \endaffil
\address La Paz Bolivia \endaddress
\email jguachallah\@gmail.com \endemail
\date October 2008 \enddate
\subjclass 46B Functional Analisis \endsubjclass
\abstract
We show $L^1$ is complemented in the dual space $L^{\infty *}$ for a
finite regular complex measure on a compact Hausdorff space
\endabstract
\endtopmatter
\document

\head Introduction \endhead

Let $X$ be a compact Hausdorff topological space with a complex regular finite measure $m$. Consider the isometric inclusion map
$$C(X)@>i>>L^\infty$$
Then, its adjoint is surjective
$$i^*:L^{\infty *}@>>> C(X)^*$$
By the Riesz representation theorem $C(X)^*$ is isometrically isomorphic to $M(X)$ the Banach space of complex regular measures on $X$.
By the Lebesgue-Radon-Nikodym theorem [R], identify the absolutly continuous measures with respect to the measure $m$, with the $L^1$ Banach space of integrable functions.
Let us then define
$$\phi:L^{\infty *}@>>> L^1$$
by
$$\phi(\lambda)=g$$
where $g$ is the Radon-Nikodym derivative of $\frac{d\mu_a}{dm}$, and $d\mu_a$ being the absolutly continous part of the measure related by the Riesz theorem to the functional $\lambda|_{C(X)}$.
\proclaim{Theorem} The short exact sequence 
$$0@>>>K(\phi)@>>>L^{\infty *}@>\phi>>L^1@>>>0$$
is split. And therefore 
$$L^{\infty *}\cong L^1 \oplus K(\phi)$$
\endproclaim
\demo{Proof} Let us call
$$\phi(\rho(g))=h$$
According to Lebesgue-Radon-Nikodym theorem, for all $f$, $\mu_a$-integrable
$$\int f d\mu_a=\int fhdm$$
Therefore, in particular for $f\in L^{\infty}$
$$\int fh dm=\int fg dm$$
And by the Hahn-Banach theorem
$$h=g\ a.e.(m)$$
\enddemo
To determine the kernel of $\phi$
\proclaim{Proposition} Kernel of $\phi$ is
$$K(\phi)=C(X)^o + M_s$$
where $C(X)^o$ is the anhilitator of $C(X)$ in $L^{\infty *}$ and $M_s$ is the singular part of a measure respect to the measure $m$.
\endproclaim
\demo{Proof} Since the kernel of the map
$$M_a\oplus M_s@>>>M_a$$
is $M_s$ and 
$$\frac{L^{\infty *}}{K(i^*)}\cong C(X)^*$$
then the kernel of $\phi$ is the sum
$$K(\phi)=C(X)^o + M_s$$
\enddemo

\end{document}